\documentclass[11pt]{amsart}
\usepackage{amsmath,amssymb,amsthm}
\usepackage[margin=1.1in]{geometry}
\usepackage{hyperref}

\newtheorem{theorem}{Theorem}
\newtheorem{lemma}[theorem]{Lemma}
\theoremstyle{definition}
\newtheorem{definition}[theorem]{Definition}
\theoremstyle{remark}
\newtheorem{remark}[theorem]{Remark}

\newcommand{\ZZ}{\mathbb{Z}}
\newcommand{\CC}{\mathbb{C}}
\newcommand{\abs}[1]{\lvert#1\rvert}
\newcommand{\norm}[1]{\lVert#1\rVert}
\newcommand{\tp}{\top}
\newcommand{\ii}{\mathrm{i}}
\DeclareMathOperator{\tr}{tr}
\DeclareMathOperator{\SL}{SL}
\DeclareMathOperator{\SU}{SU}
\DeclareMathOperator{\SO}{SO}
\DeclareMathOperator{\ord}{ord}
\newcommand{\su}{\mathfrak{su}(2)}
\newcommand{\slc}{\mathfrak{sl}_2(\CC)}

\title[A unipotent witness]{A unipotent witness: vanishing orders of word maps\\ on $\SU(2)$ and lengths in the lower central series}
\author{Daniel Keren}
\address{Department of Computer Science, University of Haifa, Haifa, Israel}
\email{dkeren@ds.haifa.ac.il}
\date{November 2026}

\begin{document}

\begin{abstract}
Let $n\ge2$, $k\ge2$ and $1\ne w\in\gamma_k(F_n)$, the $k$-th term of the lower central series of the free group
of rank $n$. We prove that a cyclically reduced conjugate of $w$ has at least $2k$ cyclic syllables, so $w$ has word length
at least $2k$; the previous linear bound, due to Malestein and Putman, is $k$. For $1\ne w\in[F_n,F_n]$ we also prove
that the word map of $w$ on $\SU(2)$ or $\SO(3)$ vanishes at the identity to order at most $T/2$, where $T$ is the
number of cyclic syllables of a cyclically reduced conjugate of $w$, with equality for $[x_1,x_2]$ and for words with six
and eight syllables. Both proofs evaluate $w$ at unipotent matrices $I+tN_i$, with $N_i$ nilpotent of rank one, and compare
the degree of the trace of the resulting matrix $M(t)\in\SL_2(\ZZ[t])$ with the order of vanishing of $M(t)-I$ at
$t=0$.
\end{abstract}

\maketitle

\section{Introduction and results}

Let $n\ge2$ and let $F_n$ be the free group on $x_1,\dots,x_n$. Write $F_n'=[F_n,F_n]$ and $\gamma_k=\gamma_k(F_n)$
for the lower central series ($\gamma_1=F_n$, $\gamma_{k+1}=[\gamma_k,F_n]$, so $\gamma_2=F_n'$), where
$[x,y]=xyx^{-1}y^{-1}$. For a word $w$, $n_i$ denotes the exponent sum of $x_i$ in $w$; $w\in F_n'$ if and only if all
$n_i=0$. A \emph{syllable} is a maximal block $x_i^{r}$, $r\ne0$, of a reduced word. The \emph{cyclic syllables} of a
cyclically reduced word are its syllables when the word is read cyclically, so that the last and the first syllable
merge if they are powers of the same generator. For example, the cyclically reduced word $x_1x_2x_1^{-1}x_2x_1$ has $5$ syllables
but $4$ cyclic syllables, since its last and first syllables merge into $x_1^2$.

We write $M_2(\CC)$ for the algebra of complex $2\times2$ matrices, $\slc\subset M_2(\CC)$ for the traceless ones,
and $\su\subset\slc$ for the Lie algebra of $\SU(2)$, the traceless skew-Hermitian matrices; $E_{pq}$ is the matrix
with $1$ in position $(p,q)$ and $0$ elsewhere, and for matrices $[X,Y]=XY-YX$. Norms of matrices are operator norms with respect
to the Euclidean norm on $\CC^2$.

\begin{definition}\label{def:ord}
For $1\ne w\in F_n$ and an integer $d\ge1$, say that $w$ \emph{vanishes to order $d$ at $I$} if there are
$C,\varepsilon_0>0$ such that
\[
\norm{w(h_1,\dots,h_n)-I}\le C\,\varepsilon^d\quad\text{whenever }0<\varepsilon\le\varepsilon_0,\ h_i\in\SU(2)\text{ and
}\norm{h_i-I}\le\varepsilon\ (1\le i\le n).
\]
This holds for $d=1$ (Lemma~\ref{lem:witness}(d) below), and if it holds for $d$ it holds for every $d'\le d$. Let
$\ord(w)\in\{1,2,\dots\}\cup\{\infty\}$ be the supremum of such $d$.
\end{definition}

By Lemma~\ref{lem:witness}(a),(d) below, $\ord(w)$ is the least $k$ for which the $k$-th homogeneous term of
$w(e^{tX_1},\dots,e^{tX_n})$ is not identically zero on $\su^n$. In particular $\ord(w)=1$ if $w\notin F_n'$, since
that term is then $\sum_in_iX_i$. The covering map $\rho\colon\SU(2)\to\SO(3)$ restricts to a diffeomorphism between
neighborhoods of $I$, so every $g\in\SO(3)$ close to $I$ is $\rho(h)$ for a unique $h\in\SU(2)$ close to $I$, and
$\norm{\rho(h)-I}$ and $\norm{h-I}$ are bounded by constant multiples of each other near $I$; since
$\rho(w(h_1,\dots,h_n))=w(\rho(h_1),\dots,\rho(h_n))$, Definition~\ref{def:ord} with $\SO(3)$ in place of $\SU(2)$
(and any norm on real $3\times3$ matrices) gives the same number.

\begin{theorem}\label{thm:A}
Let $1\ne w\in F_n'$ be a word of length $L$, and let $T$ be the number of cyclic syllables of a cyclically reduced
conjugate of $w$ (any two such conjugates are cyclic permutations of each other, so $T$ depends only on $w$). Then
\[
\ord(w)\le T/2\le L/2 .
\]
The bound is attained by $w=[x_1,x_2]$, for which $\ord(w)=2$ and $T=L=4$ in every rank $n\ge2$.
\end{theorem}

\begin{theorem}\label{thm:B}
Let $n\ge2$, $k\ge2$ and $1\ne w\in\gamma_k(F_n)$. Then a cyclically reduced conjugate of $w$ has at least $2k$
cyclic syllables; in particular the reduced word representing $w$ has length at least $2k$.
\end{theorem}

Before Theorem~\ref{thm:B}, the best lower bound for the length was $k$, due to Malestein and
Putman~\cite[Theorem~1.2]{MP} for free groups of any rank, improving the bound $k/2$ that follows from
Fox~\cite[Lemma~4.2]{Fox}; the shortest elements of $\gamma_k(F_2)$ are known to
have length at most of order $k^{1.45}$ (Elkasapy--Thom~\cite{ET}, Elkasapy~\cite{E}); see~\cite{MP2} for a recent
account of the gap. Theorem~\ref{thm:B} improves the constant of the linear lower bound from $1$ to $2$, in every
rank. Like the proof of~\cite[Theorem~1.2]{MP}, which shows that a nontrivial element of $\gamma_k$ has at least
$k$ syllables, it bounds the number of syllables and not only the length; it does not change the exponent. Its proof is purely algebraic and
does not use Section~\ref{sec:witness}. The syllable bound is attained for $k=2$, $3$ and $4$
(Remark~\ref{rem:ex}); we do not know whether the constant $2$ can be improved for large $k$. The argument uses
only the degree of $\tr M(t)$ and its order of vanishing at $t=0$; an improvement would presumably need more
information about $\tr M(t)$.

Theorem~\ref{thm:A} bounds the quantity behind the known almost laws for $\SO(3)$, nontrivial words $W$ of length
$L$ with $\norm{W(A,B)-I}\le\exp(-cL^{\delta})$ for all $A,B\in\SO(3)$: Elkasapy~\cite{E} has
$\delta=\log_2\bigl((1+\sqrt5)/2\bigr)=0.694\ldots$ and Yehuda~\cite{Y} has $\delta=0.879\ldots$, both from words whose
order at the identity grows like a power of the length. By Theorem~\ref{thm:A}, that order is at most $L/2$.

The weaker bound $\ord(w)\le L$ in Theorem~\ref{thm:A} follows, via Lemma~\ref{lem:witness}(b), from work of Glazer
and Hendel~\cite[Proposition~7.6]{GH}, who show that the symbol of a word map on a simple algebraic group, its
lowest-order nonzero term at the identity, is a Lie algebra word map of degree at most the length of the word.

Section~\ref{sec:witness} shows that $\ord(w)$ is bounded by the order of vanishing of $w(e^{tX_1},\dots,e^{tX_n})-I$
at $t=0$ for \emph{any} $X_1,\dots,X_n\in\slc$, even when the $X_i$ are not in $\su$, so that the matrices
$e^{tX_i}$ are not in $\SU(2)$. In Section~\ref{sec:unip} we take nilpotent $X_i$ of rank one, for which the word
becomes a matrix $M(t)$ over $\ZZ[t]$. The results follow from two facts about its trace. First, for cyclically reduced $w$ involving at least two generators,
the degree of $\tr M(t)$ equals the number of cyclic syllables of $w$; for $n=2$ this is a special case of the triangular
specialization of Horowitz~\cite{H}, in the form of~\cite[Lemma~4.4]{K}. Second, since $\det M=1$, the polynomial
$\tr M(t)-2$ vanishes at $t=0$ to at least twice the order of $M(t)-I$. As far as we know, these two facts have not
been combined before.

\section{Witnesses outside $\su$}\label{sec:witness}

To bound $\ord(w)$ from above it suffices to exhibit a single tuple at which $w$ is not too close to $I$. We call
$(X_1,\dots,X_n)\in\slc^n$ a \emph{witness of order $B$} for $w$, or simply a witness, if
$w(e^{tX_1},\dots,e^{tX_n})-I$ vanishes to order exactly $B<\infty$ at $t=0$; by Lemma~\ref{lem:witness}(c) below, $\ord(w)\le B$. A witness need not lie in $\su^n$: the
coefficients $c_k$ below are polynomial, so they vanish on $\su^n$ exactly when they vanish on $\slc^n$.

\begin{lemma}\label{lem:expansion}
For $X_1,\dots,X_n\in M_2(\CC)$ the map $t\mapsto\Phi_w(t)=w(e^{tX_1},\dots,e^{tX_n})$ is given by a power series in
$t$ that converges for every $t\in\CC$:
\[
\Phi_w(t)=I+\sum_{k\ge1}c_k(X_1,\dots,X_n)\,t^k ,
\]
where each $c_k:M_2(\CC)^n\to M_2(\CC)$ is a polynomial map with rational coefficients, homogeneous of degree
$k$.
\end{lemma}

\begin{proof}
Each letter contributes a factor $e^{\pm tX_i}=\sum_j(\pm tX_i)^j/j!$, a series that converges absolutely for every
$t$; hence so does their product, which has one factor for each letter of $w$. Multiplying out and collecting powers
of $t$, the coefficient of $t^k$ is a finite sum, with rational coefficients, of products of the $X_i$ of total length
$k$. The constant term is $I$.
\end{proof}

\begin{lemma}\label{lem:witness}
Let $c_k$ be the coefficient maps of Lemma~\ref{lem:expansion} for the word $w$, and let
$m=\min\{k\ge1: c_k\not\equiv0\text{ on }\su^n\}\in\mathbb N\cup\{\infty\}$ (with $m=\infty$ if the
set is empty). Then:
\begin{enumerate}
\item[(a)] $\ord(w)\le m$;
\item[(b)] $m=\min\{k\ge1: c_k\not\equiv0\text{ on }\slc^n\}$;
\item[(c)] if for some $X_1,\dots,X_n\in\slc$ one has $w(e^{tX_1},\dots,e^{tX_n})=I+t^BW(t)$ for an integer $B\ge1$
and a matrix-valued power series $W(t)$ with $W(0)\ne0$ (that is, $(X_1,\dots,X_n)$ is a witness of order $B$
for $w$), then $\ord(w)\le B$;
\item[(d)] $\ord(w)\ge m$; hence, with (a), $\ord(w)=m$.
\end{enumerate}
\end{lemma}

\begin{proof}
(a) There is nothing to prove if $m=\infty$. Otherwise choose $X_1,\dots,X_n\in\su$ with $c_m(X_1,\dots,X_n)\ne0$;
replacing each $X_i$ by $\lambda X_i$ for a real $\lambda>0$ we may assume $\norm{X_i}\le1$. For $0<\varepsilon\le1$
put $t=\varepsilon/2$. Then $e^{tX_i}\in\SU(2)$ and
\[
\norm{e^{tX_i}-I}\le t\norm{X_i}e^{t\norm{X_i}}\le(\varepsilon/2)e^{1/2}\le\varepsilon
\]
(factor $tX_i$ out of the series for $e^{tX_i}-I$), while, since $c_k\equiv0$ on $\su^n$ for $k<m$,
Lemma~\ref{lem:expansion} gives
\[
\norm{w(e^{tX_1},\dots,e^{tX_n})-I}\ \ge\ \norm{c_m}\,t^m-O(t^{m+1})\ \ge\ \tfrac12\norm{c_m}\,2^{-m}\varepsilon^m
\]
for small $\varepsilon$, with $c_m=c_m(X_1,\dots,X_n)$. So $w$ does not vanish to order $m+1$ at $I$.

(b) The matrices
\[
Z_1=\begin{pmatrix}\ii&0\\0&-\ii\end{pmatrix},\quad Z_2=\begin{pmatrix}0&1\\-1&0\end{pmatrix},\quad
Z_3=\begin{pmatrix}0&\ii\\\ii&0\end{pmatrix}
\]
lie in $\su$ and are linearly independent over $\CC$, so they form a basis of the $3$-dimensional complex space
$\slc$; every $X\in\slc$ is $X=\sum_j\xi_jZ_j$ with unique $\xi_j\in\CC$, and $X\in\su$ exactly when all $\xi_j$ are
real. In the coordinates $X_i=\sum_j\xi_{ij}Z_j$, each entry of $c_k$ is a complex polynomial in
$(\xi_{ij})\in\CC^{3n}$. If it vanishes for all real $(\xi_{ij})$, all its coefficients vanish, so it vanishes on
$\CC^{3n}$. Hence $c_k\equiv0$ on $\su^n$ implies $c_k\equiv0$ on $\slc^n$; the converse is trivial.

(c) Comparing the coefficients of $t^k$ in $I+t^BW(t)$ and in the expansion of Lemma~\ref{lem:expansion},
$c_k(X_1,\dots,X_n)=0$ for $k<B$ and $c_B(X_1,\dots,X_n)=W(0)\ne0$. So $c_B\not\equiv0$ on $\slc^n$, hence $m\le B$
by (b), and $\ord(w)\le m\le B$ by (a).

(d) Let $L$ be the length of $w$. If $h\in\SU(2)$ and $\norm{h-I}\le\varepsilon$, diagonalize $h$ by a unitary
matrix: $h$ has eigenvalues $e^{\pm\ii\theta}$ with $0\le\theta\le\pi$, and $h=e^{X}$ for some $X\in\su$ with
eigenvalues $\pm\ii\theta$. Then $\norm{X}=\theta$ and $\norm{h-I}=\abs{e^{\ii\theta}-1}=2\sin(\theta/2)\ge2\theta/\pi$,
so $\norm{X}\le\pi\varepsilon/2$. Now let $h_1,\dots,h_n\in\SU(2)$ with $\norm{h_i-I}\le\varepsilon$; by the above, $h_i=e^{X_i}$ with
$X_i\in\su$ and $\norm{X_i}\le r:=\pi\varepsilon/2$. By
Lemma~\ref{lem:expansion} with $t=1$, $w(h_1,\dots,h_n)=I+\sum_{k\ge1}c_k(X_1,\dots,X_n)$, and
$c_k(X_1,\dots,X_n)$ is a sum, over $(j_1,\dots,j_L)$ with $\sum_lj_l=k$, of products with norm at most
$\prod_lr^{j_l}/j_l!$; hence $\norm{c_k(X_1,\dots,X_n)}\le(Lr)^k/k!$. Since $c_k\equiv0$ on $\su^n$ for $k<m$,
\[
\norm{w(h_1,\dots,h_n)-I}\le\sum_{k\ge m}\frac{(Lr)^k}{k!}\le(Lr)^me^{Lr}
\le\Bigl(\frac{\pi L}{2}\Bigr)^me^{\pi L/2}\,\varepsilon^m\qquad(\varepsilon\le1).
\]
So $w$ vanishes to order $m$ at $I$. (For $m=\infty$ apply this to every finite $m'$.) Since $m\ge1$, this also
gives the case $d=1$ used in Definition~\ref{def:ord}.
\end{proof}

\section{The unipotent witness}\label{sec:unip}

Let $J=\begin{pmatrix}0&1\\-1&0\end{pmatrix}$ and $\omega(u,v)=u^{\tp}Jv=\det(u\mid v)$ for $u,v\in\ZZ^2$ (the
determinant of the matrix with columns $u,v$). Choose $u_1,\dots,u_n\in\ZZ^2$ with $\omega(u_i,u_j)\ne0$ for $i\ne j$
(for instance $u_i=(1,i)$, for which $\omega(u_i,u_j)=j-i$), and define
\[
N_i=u_iu_i^{\tp}J,\qquad g_i(t)=e^{tN_i}=I+tN_i .
\]
Since $u_i^{\tp}Ju_i=0$, we have $\tr N_i=0$ and $N_i^2=0$. So $N_i\in\slc$, the exponential series terminates after the
linear term, $g_i(t)\in\SL_2(\ZZ[t])$, and $g_i(t)^r=I+rtN_i$ for every integer $r$. For a word $w$ let
\[
M(t):=w(g_1(t),\dots,g_n(t))\in\SL_2(\ZZ[t]).
\]
Every syllable $x_i^r$ contributes the single factor $I+rtN_i$, whose entries have degree $\le1$; so the entries
of $M$ have degree at most the number of syllables of $w$.

For $n=2$, $u_1=(1,0)$, $u_2=(0,1)$ one gets $N_1=E_{12}$, $N_2=-E_{21}$. In this case the next lemma is a special
case of the triangular specialization of Horowitz~\cite{H} for words in the triangular matrices
$\begin{pmatrix}\lambda&t\\0&\lambda^{-1}\end{pmatrix}$ and $\begin{pmatrix}\mu&0\\t&\mu^{-1}\end{pmatrix}$, obtained by
setting $\lambda=\mu=1$ and replacing $t$ by $-t$ in the second matrix (see~\cite[Lemma~4.4]{K}, where the leading
coefficient is computed for general $\lambda,\mu$).

\begin{lemma}[exact degree of the trace]\label{lem:deg}
Let $w=x_{i_1}^{r_1}x_{i_2}^{r_2}\cdots x_{i_T}^{r_T}$ be a cyclically reduced word, written as its $T\ge2$ cyclic
syllables, so that $r_j\ne0$ and consecutive syllables, including the last and the first, involve different
generators ($i_j\ne i_{j+1}$, indices mod $T$). Then $\tr M(t)$ is a polynomial in $t$ of degree exactly $T$, with
leading coefficient
\[
\prod_{j=1}^Tr_j\;\prod_{j=1}^T\omega\bigl(u_{i_j},u_{i_{j+1}}\bigr)\ \ne\ 0 .
\]
\end{lemma}

\begin{proof}
$M(t)=\prod_{j=1}^T(I+r_jtN_{i_j})$. Expanding the product, a term that includes the summand $I$ from some factor
has degree at most $T-1$ in $t$; the only term of degree $T$ is $t^T\prod_jr_j\cdot N_{i_1}N_{i_2}\cdots N_{i_T}$.
Writing $v_j=u_{i_j}$ (and $v_{T+1}=v_1$),
\[
N_{i_1}\cdots N_{i_T}=v_1\,(v_1^{\tp}Jv_2)(v_2^{\tp}Jv_3)\cdots(v_{T-1}^{\tp}Jv_T)\,v_T^{\tp}J,
\]
whose trace is $\prod_{j<T}\omega(v_j,v_{j+1})\cdot v_T^{\tp}Jv_1=\prod_{j=1}^T\omega(v_j,v_{j+1})$. Each factor is
nonzero because $i_j\ne i_{j+1}$.
\end{proof}

\begin{lemma}[determinant identity]\label{lem:det}
Let $M\in\SL_2(\CC[t])$, $M=I+t^BW$ with $B\ge1$, $W\in M_2(\CC[t])$. Then
\[
\tr M-2=-\,t^{2B}\det W .
\]
\end{lemma}

\begin{proof}
For $2\times2$ matrices, $\det(I+Q)=1+\tr Q+\det Q$. With $Q=t^BW$ and $\det M=1$:
$t^B\tr W+t^{2B}\det W=0$, so $\tr M-2=\tr Q=t^B\tr W=-t^{2B}\det W$.
\end{proof}

\begin{proof}[Proof of Theorem~\ref{thm:A}]
Let $1\ne w\in F_n'$. The order is invariant under conjugation of $w$ (conjugation in $\SU(2)$ preserves
$\norm{\,\cdot-I}$), and so is the trace of $M$; so we may replace $w$ by a cyclically reduced conjugate, which
has length at most $L$, lies in $F_n'$ and is nontrivial. It involves at least two generators (a nontrivial power
of one generator has a nonzero exponent sum). A cyclic permutation of it is again a conjugate of $w$, so we may
also assume that its first and last syllables are powers of different generators. It then has the form of
Lemma~\ref{lem:deg}, with $2\le T\le L$.

By Lemma~\ref{lem:deg}, $\tr M-2$ is a nonzero polynomial of degree $T$; in particular $M\ne I$. Since $M(0)=I$,
write $M=I+t^BW$ with $B\ge1$ and $W(0)\ne0$. By Lemma~\ref{lem:det}, $t^{2B}$ divides the nonzero polynomial
$\tr M-2$, so $2B\le T$. By Lemma~\ref{lem:witness}(c), applied with $X_i=N_i$, $\ord(w)\le B\le T/2\le L/2$.

\emph{Equality for $[x_1,x_2]$.} Here $L=T=4$, so $\ord([x_1,x_2])\le2$ by what we just proved. Conversely, if
$w\in F_n'$ then $c_1(X_1,\dots,X_n)=\sum_in_iX_i=0$ because all exponent sums vanish, so $m\ge2$ and $\ord(w)\ge2$
by Lemma~\ref{lem:witness}(d).
\end{proof}

\begin{proof}[Proof of Theorem~\ref{thm:B}]
Let $\Gamma_p=\{U\in\SL_2(\ZZ[t]): U\equiv I \bmod t^p\}$. For $U=I+t^pP\in\Gamma_p$ and $V=I+t^qQ\in\Gamma_q$ (with
$P,Q\in M_2(\ZZ[t])$),
\[
UVU^{-1}V^{-1}-I=(UV-VU)\,U^{-1}V^{-1}=t^{p+q}(PQ-QP)\,U^{-1}V^{-1},
\]
and since $\Gamma_{p+q}$ is a subgroup (the kernel of reduction modulo $t^{p+q}$), $[\Gamma_p,\Gamma_q]\subset\Gamma_{p+q}$.
The homomorphism $\varphi\colon F_n\to\SL_2(\ZZ[t])$, $x_i\mapsto g_i(t)$, maps $F_n$ into $\Gamma_1$, so by induction
$\varphi(\gamma_k)=[\varphi(\gamma_{k-1}),\varphi(F_n)]\subset[\Gamma_{k-1},\Gamma_1]\subset\Gamma_k$.

Let $1\ne w\in\gamma_k$, $k\ge2$. As $\gamma_k$ is normal, we may, as in the proof of Theorem~\ref{thm:A}, replace $w$ by a cyclically reduced conjugate with
$T$ cyclic syllables, written as in Lemma~\ref{lem:deg}; it lies in $F_n'$. By Lemma~\ref{lem:deg}, $M\ne I$, so
$M=I+t^BW$ with $W(0)\ne0$, and $B\ge k$ because $M\in\Gamma_k$. As in the proof of Theorem~\ref{thm:A}, $2B\le T$.
Hence $2k\le T$, and $T$ is at most the length of the cyclically reduced conjugate, which is at most the length of
the reduced word representing $w$.
\end{proof}

\begin{remark}
Theorem~\ref{thm:B} does not use Section~\ref{sec:witness}: it uses only $\SL_2(\ZZ[t])$, the filtration
$\Gamma_1\supset\Gamma_2\supset\cdots$ and Lemmas~\ref{lem:deg}, \ref{lem:det}. The syllable bound $2k$ is attained
for $k=2$ ($[x_1,x_2]$), for $k=3$ (the word $ab^2ab^{-1}a^{-2}b^{-1}$ of Remark~\ref{rem:ex}, which lies in
$\gamma_3$ and has $6$ syllables) and for $k=4$ (the word $a^3b^2ab^{-1}a^{-3}b^{-2}a^{-1}b$ of
Remark~\ref{rem:ex}, which lies in $\gamma_4$ and has $8$ syllables).
\end{remark}

\begin{remark}[Examples, $n=2$]\label{rem:ex}
Let $a=x_1$, $b=x_2$, and take $u_1=(1,0)$, $u_2=(0,1)$, so that $N_1=E_{12}$ and $N_2=-E_{21}$. Throughout this
remark we use the witness $(N_1,-N_2)$ instead of $(N_1,N_2)$: we put $\alpha(t)=e^{tE_{12}}=g_1(t)$,
$\beta(t)=e^{tE_{21}}=g_2(-t)$ and $M(t)=w(\alpha(t),\beta(t))$. Lemmas~\ref{lem:witness}
and~\ref{lem:det} apply unchanged, and Lemma~\ref{lem:deg} holds with its leading coefficient multiplied by $(-1)^s$,
where $s$ is the number of $b$-syllables. For $w=ab^2ab^{-1}a^{-2}b^{-1}$ (found by a computer search), which lies
in $F_2'$ and has $L=8$ and $T=6$,
\[
M(t)=\begin{pmatrix}1-4t^4-4t^6&2t^3+4t^5\\-2t^3-4t^5&1+4t^4\end{pmatrix},
\]
so $\tr M-2=-4t^6$ and $M-I$ vanishes to order exactly $3$. Hence $\ord(w)\le3$.

We now show $\ord(w)\ge3$. Write $w=y_1\cdots y_L$ with $y_p\in\{a^{\pm1},b^{\pm1}\}$, let $\epsilon_p=\pm1$ be the
exponent of $y_p$, and for $X,Y\in M_2(\CC)$ let $A_p=\epsilon_pX$ if $y_p=a^{\pm1}$ and $A_p=\epsilon_pY$ if
$y_p=b^{\pm1}$. The $t^2$-coefficient of $\prod_pe^{tA_p}$ is $c_2=\frac12\sum_pA_p^2+\sum_{p<q}A_pA_q$. If
$w\in F_2'$ then $\sum_pA_p=0$, and squaring gives $\sum_pA_p^2=-\sum_{p<q}(A_pA_q+A_qA_p)$, so
$c_2=\frac12\sum_{p<q}[A_p,A_q]$. Pairs of two $a$-letters or of two $b$-letters contribute nothing, since $[X,X]=[Y,Y]=0$; and since
$\bigl(\sum_{y_p\in\{a^{\pm1}\}}\epsilon_p\bigr)\bigl(\sum_{y_q\in\{b^{\pm1}\}}\epsilon_q\bigr)=0$, the pairs with the
$b$-letter first contribute the same as those with the $a$-letter first. Hence
\[
c_2=\nu(w)\,[X,Y],\qquad \nu(w)=\sum_{\substack{p<q\\ y_p\in\{a^{\pm1}\},\ y_q\in\{b^{\pm1}\}}}\epsilon_p\epsilon_q .
\]
(For $[a,b]$, $\nu=1$; in general $\nu(w)$ is the signed area, counterclockwise positive, enclosed by the closed
lattice path of $w$, in which $a^{\pm1}$ and $b^{\pm1}$ are unit steps right/left and up/down.) For our $w$ the four
$a$-letters contribute $0,-2,1,1$, so $\nu(w)=0$. Thus $c_1\equiv c_2\equiv0$, and $\ord(w)\ge3$ by
Lemma~\ref{lem:witness}(d). So $\ord(w)=3=T/2$, and the bound $T/2$ of Theorem~\ref{thm:A} is attained.

Moreover, $w\in\gamma_3$. Indeed, $c_2$ is additive on $F_2'$, since $w(e^{tX},e^{tY})=I+O(t^2)$ there, and it
vanishes on $\gamma_3=[F_2',F_2]$: if $U=I+O(t^2)$ and $V=I+O(t)$ are power series in $t$, then, as in the proof
of Theorem~\ref{thm:B}, $UVU^{-1}V^{-1}-I=(UV-VU)\,U^{-1}V^{-1}=O(t^3)$. So $\nu$ induces a homomorphism $\gamma_2/\gamma_3\to\ZZ$ sending $[a,b]$ to $1$. In addition, $\gamma_2/\gamma_3$ is
generated by the image of $[a,b]$: $\gamma_2$ is the normal closure of $[a,b]$, since $F_2$ modulo that normal
closure is abelian, and modulo $\gamma_3$ every conjugate of $[a,b]$ equals $[a,b]$, because
$g[a,b]g^{-1}=[a,b]\,[[a,b]^{-1},g]$. So the homomorphism is an isomorphism, and $\nu(w)=0$ gives $w\in\gamma_3$. Since $w$ has $T=6$ syllables, it also shows that the bound
$T\ge2k$ of Theorem~\ref{thm:B} is attained for $k=3$.

For $k=4$, the word $w_4=a^3b^2ab^{-1}a^{-3}b^{-2}a^{-1}b$, also found by a computer search, has $L=14$ and $T=8$.
With the witness above,
\[
M(t)=\begin{pmatrix}1-6t^4-12t^6+36t^8&24t^5+36t^7\\-6t^5+12t^7&1+6t^4+12t^6\end{pmatrix},
\]
so $\tr M-2=36t^8$ and $M-I$ vanishes to order exactly $4$; hence $\ord(w_4)\le4$. Let $\xi,\eta$ be non-commuting
variables. The Magnus expansion of $w_4$, its image under $a\mapsto1+\xi$, $b\mapsto1+\eta$ (so that
$a^{-1}\mapsto1-\xi+\xi^2-\cdots$ and $b^{-1}\mapsto1-\eta+\eta^2-\cdots$) in the ring of formal power series in
$\xi,\eta$, has no terms of degree $1$, $2$ or $3$ (a direct computation, done by computer). By a theorem of
Magnus~\cite{Ma}, $\gamma_k$ consists of the elements of $F_2$ whose Magnus expansion has no terms of degree
$1,\dots,k-1$. Hence $w_4\in\gamma_4$ (and $w_4\notin\gamma_5$ by Theorem~\ref{thm:B}, since $T<10$), so the bound
$T\ge2k$ of Theorem~\ref{thm:B} is attained for $k=4$ as well. Finally, if $U=I+O(t^p)$ and $V=I+O(t^q)$ are power
series in $t$, then $UVU^{-1}V^{-1}=I+O(t^{p+q})$, as above; by induction on $k$, every $u\in\gamma_k$ satisfies
$u(e^{tX},e^{tY})=I+O(t^k)$. For $w_4\in\gamma_4$ this gives $c_1\equiv c_2\equiv c_3\equiv0$, so
$\ord(w_4)\ge4$ by Lemma~\ref{lem:witness}(d), and $\ord(w_4)=4=T/2$.
\end{remark}

\subsection*{Acknowledgments}
I thank Gal Yehuda for introducing me to almost laws in $\SU(2)$ and for helpful discussions.

The AI assistant Claude (Anthropic; model Opus 5.5) was used to find and analyze the examples of Remark~\ref{rem:ex}. For order $5$, no words with $T=10$ were found 
(the search covered words of length at most $14$ in two generators).

\end{document}